\documentclass[12pt,dvips]{amsart}
\usepackage{bbm}
\usepackage{epic}
\usepackage{eepic}
\usepackage{graphicx}
\usepackage{url}

\newtheoremstyle{Italic}{\baselineskip}{\baselineskip}{\normalfont\itshape}%
   {0pt}{\scshape}{. }{ }%
   {\thmnumber{#2. }\thmname{#1}\thmnote{ \normalfont #3}}
\newtheoremstyle{swapItalic}{.2\baselineskip}{.2\baselineskip}%
   {\normalfont\itshape}{0pt}{\normalfont\scshape}{ }{ }%
   {\thmname{#1}\thmnumber{ (#2) }\thmnote{ \normalfont #3}}
\newtheoremstyle{Roman}{\baselineskip}{\baselineskip}{\normalfont}%
   {0pt}{\scshape}{. }{ }%
   {\thmnumber{#2. }\thmname{#1}\thmnote{ \normalfont #3}}
\theoremstyle{Italic}
\newtheorem{Theorem}{Theorem}[section]
\newtheorem{Corollary}[Theorem]{Corollary}

\newtheorem{Proposition}[Theorem]{Proposition}
\theoremstyle{swapItalic}

\theoremstyle{Roman}
\newtheorem{Definition}[Theorem]{Definition}

\newcommand{\vect}[1]{\mathbf{#1}} 
\newcommand{\nv}{{\mathbf 0}} 
 
\newcommand{\eins}{{\mathbbm 1}}

\newcommand{\Z}{{\mathbb Z}}
\newcommand{\R}{{\mathbb R}}
\newcommand{\suchthat}{\thinspace : \thinspace}
\newcommand{\aff}{\operatorname{aff}} 
\newcommand{\conv}{\operatorname{conv}} 
 
\newcommand{\vol}{\operatorname{vol}}

\newcommand{\ehr}{\mbox{\rm\bf Ehr}} 
\newcommand{\ehrint}{\ehr^\circ} 
\newcommand{\EHR}{\mathfrak{Ehr}}
\newcommand{\sx}{\mathbf{S}}
\newcommand{\refldim}{\operatorname{refldim}}
\newcommand{\join}{*}

\begin{document}

\title[The reflexive dimension of a lattice polytope]{The reflexive
  dimension of a lattice\makebox[0pt]{\raisebox{1in}{\normalfont
      DUKE-CGTP-04-07}} polytope}

\author{Christian 
  Haase}
\thanks{The first author was supported in part by NSF grant
  DMS--0200740.}
\address{Duke Math \\ Durham, NC 27708-0320}
\email{haase@math.duke.edu}

\author{Ilarion V. Melnikov}
\thanks{The second author was supported in part by NSF grants
  DMS-0074072 and DMS-0301476.}
\address{Duke Physics \\ Durham, NC 27708-0305}
\email{lmel@phy.duke.edu}

\date{June 23, 2004}

\begin{abstract}
  The reflexive dimension $\refldim(P)$ of a lattice polytope $P$ is
  the minimal $d$ so that $P$ is the face of some $d$-dimensional
  reflexive polytope.

  We show that $\refldim(P)$ is finite for every $P$, and give bounds
  for $\refldim(kP)$ in terms of $\refldim(P)$ and $k$.
\end{abstract}

\maketitle
\section{Introduction}
Reflexive polytopes were originally defined with theoretical physics
applications in mind.
In string theory, reflexive polytopes and the associated toric
varieties play a crucial role in the most quantitatively predictive
form of the mirror symmetry conjecture \cite{BatyrevMirror,
  MorrisonPlesser}.
Aside from such physical uses, we want to advertise their study as
interesting combinatorial objects. They enjoy a variety of interesting
combinatorial properties~\cite{Haase24,Nill} which are not well
understood. In this note we explore how restrictive the condition of
reflexivity turns out to be. We introduce the notion of reflexive
dimension of an arbitrary lattice polytope which could be a starting
point to study the questions like
\begin{itemize}
\item Which polytopes are reflexive?
\item Can we get a handle on possible or impossible combinatorial
  types?
\item Can we find upper or lower bounds for the $f$-vector or Ehrhart
  coefficients?
\end{itemize}
A lattice polytope $P$ is the convex hull in $\R^d$ of a finite
set of lattice points, i.e., points in $\Z^d$.  Its dimension
$\dim(P)$ is the dimension of its affine span $\aff(P)$ as an
affine space.  We will identify lattice equivalent lattice polytopes,
where two lattice polytopes $P$ and $P'$ are {\em lattice equivalent}
if there exists an affine map $\aff(P) \to \aff(P')$ that maps $\Z^d
\cap \aff(P)$ bijectively onto $\Z^{d'} \cap \aff(P')$, and which maps
$P$ to $P'$.  Every lattice polytope is lattice equivalent to a
full-dimensional one, and a full-dimensional lattice polytope has a
unique presentation
\begin{equation*}
  P = \{ \vect{x} \in \R^d \suchthat \langle \vect{y}_i , \vect{x}
  \rangle \ge c_i \text{ for } i = 1,\ldots,k \} ,
\end{equation*}
where the $\vect{y}_i$ are primitive elements of the dual lattice
$(\Z^d)^\vee$, the $c_i$ are integers, and $k$ is minimal. This 
system of inequalities will also be referred to as 
$\vect{A} \vect{x} \ge \vect{c}.$
\begin{Definition}
  A lattice polytope $P = \{ \vect{x} \suchthat \vect{A} \vect{x} \ge
  \vect{c} \}$ with interior lattice point $\vect{x}_0$ is {\em
    reflexive\/} if $\vect{Ax}_0 - \vect{c} = \eins$, where $\eins$ is
  the all-one vector $(1,\ldots,1)^t$.
\end{Definition}
It follows that reflexive polytopes have precisely one interior
lattice point which lies in an adjacent lattice hyperplane to any
facet. This is sometimes described as ``all facets are distance one
from the interior lattice point.''

\begin{figure}[htbp]
\label{fig:reflexive2d}
\begin{center}
\setlength{\unitlength}{0.00016667in} 
\begingroup\makeatletter\ifx\SetFigFont\undefined%
\gdef\SetFigFont#1#2#3#4#5{%
  \reset@font\fontsize{#1}{#2pt}%
  \fontfamily{#3}\fontseries{#4}\fontshape{#5}%
  \selectfont}%
\fi\endgroup%
{\renewcommand{\dashlinestretch}{30}
\begin{picture}(20566,9781)(0,-10)
\drawline(1283,9683)(2483,8483)(83,7283)(1283,9683)
\drawline(4883,9683)(3683,7283)(6083,7283)(4883,9683)
\drawline(8483,9683)(7283,8483)(8483,7283)
        (9683,7283)(8483,9683)
\drawline(4883,6083)(3683,4883)(4883,3683)
        (6083,3683)(6083,4883)(4883,6083)
\drawline(7283,6083)(7283,3683)(8483,3683)
        (9683,4883)(9683,6083)(7283,6083)
\drawline(12083,6083)(10883,3683)(13283,3683)
        (13283,6083)(12083,6083)
\drawline(14483,3683)(14483,4883)(15683,6083)
        (16883,4883)(16883,3683)(14483,3683)
\drawline(83,6083)(83,4883)(1283,3683)
        (2483,3683)(2483,4883)(1283,6083)(83,6083)
\drawline(83,3683)(83,83)(3683,83)(83,3683)
\drawline(6083,2483)(4883,1283)(4883,83)
        (8483,83)(6083,2483)
\drawline(10883,2483)(9683,83)(13283,83)(10883,2483)
\drawline(15683,2483)(14483,83)(16883,83)
        (16883,1283)(15683,2483)
\drawline(10883,9683)(12083,9683)(14483,7283)
        (10883,7283)(10883,9683)
\drawline(15683,7283)(18083,9683)(20483,7283)(15683,7283)
\drawline(19283,2483)(18083,1283)(19283,83)
        (20483,1283)(19283,2483)
\drawline(18083,6083)(18083,3683)(20483,3683)
        (20483,6083)(18083,6083)
\put(20483,7283){\circle*{150}}
\put(19283,7283){\circle*{150}}
\put(18083,7283){\circle*{150}}
\put(16883,7283){\circle*{150}}
\put(15683,7283){\circle*{150}}
\put(20483,8483){\circle*{150}}
\put(19283,8483){\circle*{150}}
\put(18083,8483){\circle*{150}}
\put(16883,8483){\circle*{150}}
\put(15683,8483){\circle*{150}}
\put(20483,9683){\circle*{150}}
\put(19283,9683){\circle*{150}}
\put(18083,9683){\circle*{150}}
\put(16883,9683){\circle*{150}}
\put(15683,9683){\circle*{150}}
\put(20483,3683){\circle*{150}}
\put(19283,3683){\circle*{150}}
\put(18083,3683){\circle*{150}}
\put(16883,3683){\circle*{150}}
\put(15683,3683){\circle*{150}}
\put(20483,4883){\circle*{150}}
\put(19283,4883){\circle*{150}}
\put(18083,4883){\circle*{150}}
\put(16883,4883){\circle*{150}}
\put(15683,4883){\circle*{150}}
\put(20483,6083){\circle*{150}}
\put(19283,6083){\circle*{150}}
\put(18083,6083){\circle*{150}}
\put(16883,6083){\circle*{150}}
\put(15683,6083){\circle*{150}}
\put(20483,83){\circle*{150}}
\put(19283,83){\circle*{150}}
\put(18083,83){\circle*{150}}
\put(16883,83){\circle*{150}}
\put(15683,83){\circle*{150}}
\put(20483,1283){\circle*{150}}
\put(19283,1283){\circle*{150}}
\put(18083,1283){\circle*{150}}
\put(16883,1283){\circle*{150}}
\put(15683,1283){\circle*{150}}
\put(20483,2483){\circle*{150}}
\put(19283,2483){\circle*{150}}
\put(18083,2483){\circle*{150}}
\put(16883,2483){\circle*{150}}
\put(15683,2483){\circle*{150}}
\put(83,9683){\circle*{150}}
\put(83,8483){\circle*{150}}
\put(83,7283){\circle*{150}}
\put(83,6083){\circle*{150}}
\put(83,4883){\circle*{150}}
\put(83,3683){\circle*{150}}
\put(83,2483){\circle*{150}}
\put(83,1283){\circle*{150}}
\put(1283,9683){\circle*{150}}
\put(1283,8483){\circle*{150}}
\put(1283,7283){\circle*{150}}
\put(1283,6083){\circle*{150}}
\put(1283,4883){\circle*{150}}
\put(1283,2483){\circle*{150}}
\put(1283,1283){\circle*{150}}
\put(2483,9683){\circle*{150}}
\put(2483,8483){\circle*{150}}
\put(2483,7283){\circle*{150}}
\put(2483,6083){\circle*{150}}
\put(2483,4883){\circle*{150}}
\put(2483,3683){\circle*{150}}
\put(2483,2483){\circle*{150}}
\put(2483,1283){\circle*{150}}
\put(3683,9683){\circle*{150}}
\put(3683,8483){\circle*{150}}
\put(3683,7283){\circle*{150}}
\put(3683,6083){\circle*{150}}
\put(3683,4883){\circle*{150}}
\put(3683,3683){\circle*{150}}
\put(3683,2483){\circle*{150}}
\put(4883,9683){\circle*{150}}
\put(4883,8483){\circle*{150}}
\put(4883,7283){\circle*{150}}
\put(4883,6083){\circle*{150}}
\put(4883,4883){\circle*{150}}
\put(4883,3683){\circle*{150}}
\put(4883,2483){\circle*{150}}
\put(4883,1283){\circle*{150}}
\put(6083,9683){\circle*{150}}
\put(6083,8483){\circle*{150}}
\put(6083,7283){\circle*{150}}
\put(6083,6083){\circle*{150}}
\put(6083,4883){\circle*{150}}
\put(6083,3683){\circle*{150}}
\put(6083,2483){\circle*{150}}
\put(6083,1283){\circle*{150}}
\put(7283,9683){\circle*{150}}
\put(7283,8483){\circle*{150}}
\put(7283,7283){\circle*{150}}
\put(7283,6083){\circle*{150}}
\put(7283,4883){\circle*{150}}
\put(7283,3683){\circle*{150}}
\put(7283,2483){\circle*{150}}
\put(7283,1283){\circle*{150}}
\put(8483,9683){\circle*{150}}
\put(8483,8483){\circle*{150}}
\put(8483,7283){\circle*{150}}
\put(8483,6083){\circle*{150}}
\put(8483,4883){\circle*{150}}
\put(8483,3683){\circle*{150}}
\put(8483,2483){\circle*{150}}
\put(8483,1283){\circle*{150}}
\put(9683,9683){\circle*{150}}
\put(9683,8483){\circle*{150}}
\put(9683,7283){\circle*{150}}
\put(9683,6083){\circle*{150}}
\put(9683,4883){\circle*{150}}
\put(9683,3683){\circle*{150}}
\put(9683,2483){\circle*{150}}
\put(9683,1283){\circle*{150}}
\put(10883,9683){\circle*{150}}
\put(10883,8483){\circle*{150}}
\put(10883,7283){\circle*{150}}
\put(10883,6083){\circle*{150}}
\put(10883,4883){\circle*{150}}
\put(10883,3683){\circle*{150}}
\put(10883,2483){\circle*{150}}
\put(10883,1283){\circle*{150}}
\put(12083,9683){\circle*{150}}
\put(12083,8483){\circle*{150}}
\put(12083,7283){\circle*{150}}
\put(12083,6083){\circle*{150}}
\put(12083,4883){\circle*{150}}
\put(12083,3683){\circle*{150}}
\put(12083,2483){\circle*{150}}
\put(12083,1283){\circle*{150}}
\put(13283,9683){\circle*{150}}
\put(13283,8483){\circle*{150}}
\put(13283,7283){\circle*{150}}
\put(13283,6083){\circle*{150}}
\put(13283,4883){\circle*{150}}
\put(13283,3683){\circle*{150}}
\put(13283,1283){\circle*{150}}
\put(14483,9683){\circle*{150}}
\put(14483,8483){\circle*{150}}
\put(14483,7283){\circle*{150}}
\put(14483,6083){\circle*{150}}
\put(14483,4883){\circle*{150}}
\put(14483,3683){\circle*{150}}
\put(14483,2483){\circle*{150}}
\put(14483,1283){\circle*{150}}
\put(18083,83){\circle*{150}}
\put(16883,83){\circle*{150}}
\put(15683,83){\circle*{150}}
\put(14483,83){\circle*{150}}
\put(13283,83){\circle*{150}}
\put(12083,83){\circle*{150}}
\put(10883,83){\circle*{150}}
\put(9683,83){\circle*{150}}
\put(8483,83){\circle*{150}}
\put(7283,83){\circle*{150}}
\put(6083,83){\circle*{150}}
\put(4883,83){\circle*{150}}
\put(3683,83){\circle*{150}}
\put(2483,83){\circle*{150}}
\put(1283,83){\circle*{150}}
\put(83,83){\circle*{150}}
\put(1283,3683){\circle*{150}}
\put(13283,2483){\circle*{150}}
\put(3683,1283){\circle*{150}}
\end{picture}
}
\caption{All reflexive polygons (up to equivalence).}
\end{center}
\end{figure}
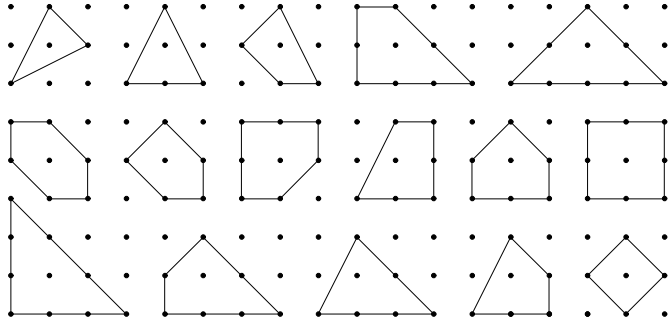

The existence of a unique interior lattice point implies that there is
only a finite number of equivalence classes of reflexive polytopes in
any given dimension~\cite{lagariasZiegler}. The only one dimensional
reflexive polytope is a segment of length $2$. In two dimensions there
are the $16$ reflexive polygons given in Figure~\ref{fig:reflexive2d}.
Three and four dimensional reflexive polytopes have been classified by
Maximilian Kreuzer and Harald Skarke
\cite{KreuzerSkarke3d,KreuzerSkarke4d}. There are $4,319$ in dimension
$3$ and  $473,800,776$ in dimension $4$.

\section{Reflexive polytopes and reflexive dimension}
The condition that a polytope be reflexive has some rather
remarkable consequences, as we now discuss. We begin with
some definitions.
\begin{itemize}
\item If $P$ is a full-dimensional polytope (not necessarily lattice)
  with the origin $\nv$ in the interior, then the polar dual
  \begin{equation*}
    P^\vee = \{ \vect{y} \in (\R^d)^* \suchthat \langle \vect{y} ,
    \vect{x} \rangle \ge -1 \text{ for all } \vect{x} \in P \}
  \end{equation*}
  is again a full-dimensional polytope with $\nv$ in the
  interior (compare~\cite[Section IV.1]{BarvinokBook}).
\item The volume $\vol(P)$ of a lattice polytope is always normalized
  with respect to the unimodular ($\dim P$)--simplex in $\aff(P) \cap
  \Z^d$.
\item The Ehrhart polynomial of a lattice polytope $P$, counts lattice
  points in dilations of $P$ (compare~\cite[Section
  VIII.5]{BarvinokBook}). \\
  $
  \ehr(P,k) := |k\, P \cap \Z^d|$ , $
    \ehrint(P,k) := |( k\, P)^\circ \cap \Z^d|$ for non-negative
    integers $k$. 
  These are, in fact, polynomials so that we can evaluate for negative
  $k$:
  $
  \ehrint(P,k) = (-1)^{\dim P} \ehr(P,-k)
  $
  . Because $\ehr$ is a polynomial, its generating function can be
  written as a rational function
  \begin{equation*}
    \sum_{k\ge0} \ehr(P,k) \; t^k \ =: \ \frac{\EHR(P,t)}{(1-t)^{\dim
        P+1}}
  \end{equation*}
\item A lattice polytope $P$ defines an ample line bundle
  $L_P$ on a projective toric variety $X_P$. (See, e.g.,
  \cite[Section~3.4]{Fulton}.)
  If $L_P$ is very ample, it provides an
  embedding $X_P \hookrightarrow \mathbb{P}^{r-1}$, where $r = |P
  \cap \Z^d|$. So we can think of $X_P$ canonically sitting in
  projective space.
\end{itemize}
The reader is invited to add some more equivalences to the
following.
\begin{Theorem}[\cite{BatyrevMirror,HibiReflexive}]
  \label{thm:refpol}
  Let $P$ be a full--dimensional lattice polytope with unique
  interior lattice point $\nv$. Then the following conditions are 
  equivalent:
  \renewcommand{\theenumi}{\roman{enumi}}
  \begin{enumerate}
    \item P is reflexive.
    \item The polar dual $P^\vee$ is a lattice polytope.
    \item $\vol(P) = \sum \vol(F)$, the sum ranging over all
      facets (codimension one faces) $F$ of $P$. \label{thm:refpol:vol}
    \item $\ehr(P,k) = \ehrint(P,k+1)$ for all $k$.
    \item $\EHR(P,1/t) = (-1)^{d+1} t \EHR(P,t)$
    \item The projective toric variety $X_P$ defined by $P$ is Fano.
    \item Every generic hyperplane section of $X_P$ is
      Calabi--Yau.\footnote{If the line bundle $L_P$ is not very
        ample, one can still define ``$P$-generic hypersurfaces'' of
        $X_P$ to generalize the notion of generic hyperplane section.}
  \end{enumerate}
\end{Theorem}
These conditions are not as restrictive as one might expect based on
the fact that in each dimension there is a finite number of reflexive
polytopes. More precisely, we have the following proposition.
\begin{Proposition}\label{prop:1}
  Every lattice polytope is lattice equivalent to a face of 
  some reflexive polytope.
\end{Proposition}
Proposition~\ref{prop:1} motivates us to define the reflexive 
dimension of a lattice polytope.
\begin{Definition}\label{def:refldim}
  Let $P$ be a lattice polytope. Then its {\em reflexive dimension\/}
  is the smallest $d$ such that $P$ is lattice equivalent to a face of 
  a reflexive $d$--polytope.
\end{Definition}
\begin{proof}[Proof of Proposition~\ref{prop:1}]
  Let $P = \{ \vect{x} \in \R^d : \vect{A}\vect{x} \ge \vect{c} \}$ be
  a lattice polytope defined by $k$ inequalities. Suppose that
  $\nv$ is an interior point of $P$.  
  Then, $c_i \le -1$, $i = 1,\ldots, k$. If equality holds then
  $P$ is already reflexive. Otherwise, we will construct a
  $(d+1)$--polytope $P'$ with strictly bigger $\vect{c}'$.

  Suppose $c_k < -1$, and introduce one more variable $x_{d+1}$.
  Then consider the polytope
  \begin{align*}
    P' = \{ (\vect{x},x_{d+1}) \in \R^{d+1} \suchthat & \langle
    \vect{y}_i , \vect{x} \rangle \ge c_i \text{ for } i =
    1,\ldots,k-1 \text{, and } \\
    & \langle \vect{y}_k , \vect{x} \rangle - x_{d+1} \ge c_k + 1
    \text{, and } \\
    & x_{d+1} \ge -1 \} .
  \end{align*}
\begin{figure}[htbp]
\label{fig:wedge}
  \begin{center}

\setlength{\unitlength}{0.00033333in}
\begingroup\makeatletter\ifx\SetFigFont\undefined%
\gdef\SetFigFont#1#2#3#4#5{%
  \reset@font\fontsize{#1}{#2pt}%
  \fontfamily{#3}\fontseries{#4}\fontshape{#5}%
  \selectfont}%
\fi\endgroup%
{\renewcommand{\dashlinestretch}{30}
\begin{picture}(3308,3942)(0,-10)
\put(83,244){\blacken\ellipse{150}{150}}
\put(83,244){\ellipse{150}{150}}
\put(2483,244){\blacken\ellipse{150}{150}}
\put(2483,244){\ellipse{150}{150}}
\put(683,244){\blacken\ellipse{150}{150}}
\put(683,244){\ellipse{150}{150}}
\put(1883,244){\blacken\ellipse{150}{150}}
\put(1883,244){\ellipse{150}{150}}
\put(83,1444){\blacken\ellipse{150}{150}}
\put(83,1444){\ellipse{150}{150}}
\put(683,1444){\blacken\ellipse{150}{150}}
\put(683,1444){\ellipse{150}{150}}
\put(1283,1444){\blacken\ellipse{150}{150}}
\put(1283,1444){\ellipse{150}{150}}
\put(1883,1444){\blacken\ellipse{150}{150}}
\put(1883,1444){\ellipse{150}{150}}
\put(2483,1444){\blacken\ellipse{150}{150}}
\put(2483,1444){\ellipse{150}{150}}
\put(83,2044){\blacken\ellipse{150}{150}}
\put(83,2044){\ellipse{150}{150}}
\put(683,2044){\blacken\ellipse{150}{150}}
\put(683,2044){\ellipse{150}{150}}
\put(1283,2044){\whiten\ellipse{150}{150}}
\put(1283,2044){\ellipse{150}{150}}
\put(1883,2044){\blacken\ellipse{150}{150}}
\put(1883,2044){\ellipse{150}{150}}
\put(2483,2044){\blacken\ellipse{150}{150}}
\put(2483,2044){\ellipse{150}{150}}
\put(83,2644){\blacken\ellipse{150}{150}}
\put(83,2644){\ellipse{150}{150}}
\put(683,2644){\blacken\ellipse{150}{150}}
\put(683,2644){\ellipse{150}{150}}
\put(1283,2644){\blacken\ellipse{150}{150}}
\put(1283,2644){\ellipse{150}{150}}
\put(1883,2644){\blacken\ellipse{150}{150}}
\put(1883,2644){\ellipse{150}{150}}
\put(2483,2644){\blacken\ellipse{150}{150}}
\put(2483,2644){\ellipse{150}{150}}
\put(83,3244){\blacken\ellipse{150}{150}}
\put(83,3244){\ellipse{150}{150}}
\put(683,3244){\blacken\ellipse{150}{150}}
\put(683,3244){\ellipse{150}{150}}
\put(1283,3244){\blacken\ellipse{150}{150}}
\put(1283,3244){\ellipse{150}{150}}
\put(1883,3244){\blacken\ellipse{150}{150}}
\put(1883,3244){\ellipse{150}{150}}
\put(2483,3244){\blacken\ellipse{150}{150}}
\put(2483,3244){\ellipse{150}{150}}
\put(83,3844){\blacken\ellipse{150}{150}}
\put(83,3844){\ellipse{150}{150}}
\put(683,3844){\blacken\ellipse{150}{150}}
\put(683,3844){\ellipse{150}{150}}
\put(1283,3844){\blacken\ellipse{150}{150}}
\put(1283,3844){\ellipse{150}{150}}
\put(1883,3844){\blacken\ellipse{150}{150}}
\put(1883,3844){\ellipse{150}{150}}
\put(2483,3844){\blacken\ellipse{150}{150}}
\put(2483,3844){\ellipse{150}{150}}
\path(83,244)(2483,244)
\path(83,1444)(2483,1444)
\path(83,1444)(83,3844)(2483,1444)
\dottedline{90}(2483,1444)(2483,3844)
\put(3308,94){\makebox(0,0)[lb]{\smash{{{\SetFigFont{10}{12.0}{\rmdefault}{\mddefault}{\updefault}$P$}}}}}
\put(833,2194){\makebox(0,0)[lb]{\smash{{{\SetFigFont{10}{12.0}{\rmdefault}{\mddefault}{\updefault}$P'$}}}}}
\put(1283,244){\whiten\ellipse{150}{150}}
\put(1283,244){\ellipse{150}{150}}
\end{picture}
}

\caption{$P$ and $P'$}
  \end{center}
\end{figure}
  Combinatorially, as seen in Figure \ref{fig:wedge}, $P'$ is the 
  wedge of $P$ over the facet 
  $\langle \vect{y}_m , \vect{x} \rangle = c_k$. 
  It has $P \times \{-1\}$ as a facet.
  Iterating this construction, we add $\| \eins + \vect{c} \|_1$
  dimensions to finally obtain a reflexive polytope. 
  
  If $P$ does not have an interior point, we can first embed $P$ as a
  face of a lattice polytope $P'$ such that $P'$ has an interior
  point and then apply the procedure above to $P'$.
\end{proof}

\section{The reflexive dimension of a line segment}
What is the reflexive dimension of a given polytope $P$?
To answer this question it is natural to start by determining the 
reflexive dimension of some sample polytopes. While the $0$-$1$-cube 
and the standard simplex are facets of reflexive polytopes, the 
question quickly gets subtle, when it comes to, e.g., multiples of 
these. In this section we will give bounds for the simplest of all 
polytopes, the segment of length $\ell$.
\subsection{Small \bfseries \itshape d}
The unique reflexive one dimensional polytope has edge length $2$, 
and by inspection of Figure \ref{fig:reflexive2d}, it is clear that
the edge lengths realized in $2$ dimensions are $\{1,2,3,4\}$. 
A search of the lists of three and four dimensional reflexive 
polytopes has yielded the following edge lengths.\footnote{The 
data is available at \url{http://tph16.tuwien.ac.at/~kreuzer/CY/}. 
The four-dimensional list can be searched by using the program
PALP~\cite{PALP}.}
\begin{Proposition}
  There is a reflexive three-polytope which has an edge of length
  $\ell$ if and only if $\ell \in \{1, \ldots, 10, 12\}$.

  There is a reflexive four-polytope which has an edge of length
  $\ell$ if and only if $\ell \in \{1, \ldots, 54, 56, 57, 58,
  60, 63, 64, 66, 70, 72, 78, 84\}$.
\end{Proposition}
So in particular, $\refldim([0,11]) = 4 > 3 = \refldim([0,12])$, and
the reflexive dimension of line segments is not monotone in $\ell$, as
could also be seen from $\refldim([0,1]) = 2 > 1 = \refldim([0,2])$.
There is precisely one reflexive three-polytope with an edge of length
$12$, and precisely one reflexive four-polytope with an edge of length
$84$.
\subsection{Large \bfseries \itshape d}
Having examined the low dimensional case, we would now like to find a 
bound on $\refldim(\left[0,\ell\right])$ for large $\ell$.

In order to obtain an asymptotic lower bound, we use a bound for the
volume of a lattice polytope which contains exactly one interior
lattice point.  
\begin{Theorem}[\cite{lagariasZiegler,Oleg}]
  If $Q$ is a $D$-polytope which contains exactly one interior
  lattice point, then
  \begin{equation*}
    \vol(Q) \le 14^{ \left(2^{D+1} D\right)} \cdot D\,!
  \end{equation*}
\end{Theorem}

If the $d$-polytope $P$ is a face of the reflexive $D$--polytope
$Q$, then $\vol(P) \le \vol(Q)$, and using the above bound we
immediately obtain a (crude) lower bound:
\begin{Corollary}
  There is a universal constant $m$ such that the reflexive dimension
  of the segment of length $\ell$ is at least $m \log \log \ell$.
\end{Corollary}

To find an upper bound,  we need to construct reflexive
polytopes with long edges. For a sequence $\vect{a} = (a_1, \ldots,
a_d)$ of positive integers, consider the simplex $\sx(\vect{a})$ which
is the convex hull of the origin and multiples $a_i \vect{e}_i$ of the
standard unit vectors.  Particularly efficient are the simplices of 
Micha Perles, J\"org Wills, and Joseph Zaks~\cite{pwz}. Let $t_1 = 2$ 
and $t_{i+1} = t_i^2 - t_i + 1$.\footnote{The sequence \url{A000058} 
in the OEIS~\cite{OEIS}} 
Then the simplex $\sx(t_1, \ldots,t_d)$ is reflexive with
interior lattice point $\eins$. It has a segment of length $t_d >
2^{2^{d-2}}$ as a face. In fact, the slight modification $\sx(t_1,
\ldots,t_{d-1}, 2t_d-2)$ is also reflexive.  Note that in 
dimensions $1,2,3,4$, these modified simplices are the 
unique polytopes that realize the longest edge length.

This is not enough to establish an upper bound $M \log \log \ell$,
since, as seen from the computational results above,  the 
reflexive dimension is not monotone in $\ell$.  The method 
presented in the
proof of Proposition~\ref{prop:1} provides the upper 
bound $\ell - 1$: it realizes the segment as a face of the $\ell$ 
times dilated standard $(\ell - 1)$--simplex. A better bound is
obtained as follows.
\begin{Proposition} \label{prop:upper}
  There is a universal constant $M$ such that the reflexive dimension
  of the segment of length $\ell$ is at most $M \sqrt{ \log \ell }$.
\end{Proposition}  
Luckily for us, Michael Vose has provided a theorem that reduces the
proof to a standard argument.
\begin{Theorem}[{\cite[Theorem 1]{Vose}}]
  There exists an increasing sequence $N_k$ of positive integers such
  that any integer $1<m<N_k$ is the sum of not more than $O(\sqrt{\log
    N_{k-1}})$ distinct divisors of $N_k$.
\end{Theorem}
\begin{proof}[Proof of Proposition~\ref{prop:upper}]
  Given $\ell$, choose $k$ so that $N_{k-1} < \ell \le N_k$. Now,
  write
  \begin{equation*}
    N_k \ell - N_k - 1 = m \ell + n
  \end{equation*}
  with $0 \le n < \ell$. Then both $m$ and $n$ are smaller than $N_k$,
  and we can write
  \begin{equation*}
    m = \sum_{i=1}^r m_i \quad \text{ and } \quad n = \sum_{j=1}^s n_j
    ,
  \end{equation*}
  where the $m_i$ and the $n_j$ divide $N_k$, and $r+s = 
  O(\sqrt{\log \ell \,})$. Using these decompositions, we can define
  integers $a_i=N_k/m_i$ and $b_j=N_k\ell/n_j$. The claim is, that the
  simplex $\sx(\vect{a}, \vect{b}, \ell)$ is reflexive with interior
  point $\eins$. (It is an $(r+s+1)$-dimensional simplex.)
  
  The only facet in question is the one which does not contain the
  origin. The integral functional $\vect{f} = \sum_{i=1}^r m_i \ell
  \vect{e}^*_i + \sum_{j=1}^s n_j \vect{e}^*_{r+j} + N_k
  \vect{e}^*_{r+s+1}$ evaluates to $N_k \ell$ on the vertices $a_i
  \vect{e}_i$, on $b_j \vect{e}_{r+j}$, and on $\ell
  \vect{e}_{r+s+1}$.
  That same functional evaluates on $\eins$ to
  \begin{equation*}
    \sum_{i=1}^r m_i \ell + \sum_{j=1}^s n_j + N_k = m \ell + n + N_k
    = N_k \ell - 1 .
  \end{equation*}
\end{proof}
It is conceivable that both bounds are sharp, in the sense that we do
know of a sequence of lengths whose reflexive dimension behaves like
$\log \log \ell$, and there might be a different sequence of lengths
whose reflexive dimension behaves like $\sqrt{\log \ell \,}$.
\section{Products, Joins, Dilations}
In this section, we collect some more and some less trivial
observations how the reflexive dimension behaves with respect to
standard operations on polytopes.

If $P \subset \R^d$ and $P' \subset \R^{d'}$ are polytopes, their join
$P \join P'$ is the convex hull in $\R^{d+d'+1}$ of $P \times \{\nv\}
\times 0$ and $\{\nv\} \times P' \times 1$.
\begin{Proposition}
  \begin{align*}
    \refldim(P \times P') &\le \refldim(P) + \refldim(P') \\
    \refldim(P \join P') &\le \refldim(P) + \refldim(P')
  \end{align*}
\end{Proposition}
\begin{proof}
  If $P$, $P'$ are faces of reflexive polytopes $Q$, $Q'$
  respectively, then $P \times P'$ is a face of $Q \times Q'$, and
  (unless $P=Q$ or $P=Q'$) $P \join P'$ is a face of $\conv(Q \times
  \{\nv\} \cup \{\nv\} \times Q')$.
\end{proof}
We can employ the same method we used for $[0,\ell] = \ell \, [0,1]$
to deal with dilations of polytopes in general.
\begin{Proposition}\label{prop:polymult}
  Suppose that the simplex $\sx(a_1, \ldots, a_r)$ is reflexive, and
  that $P$ is a reflexive $s$-polytope with interior point $\nv$. Then
  the convex hull of $(a_1 \vect{e}_1,\nv), \ldots, (a_{r-1}
  \vect{e}_{r-1},\nv)$, and $(a_r \vect{e}_r, a_r P)$ is a reflexive
  $(r+s)$-polytope which contains a face (equivalent to) $a_r P$.
\end{Proposition}
\begin{proof}
  Suppose $P$ is given by $\{ \vect{y} \suchthat \vect{Ay} \ge -\eins
  \}$.
  Let $Q_1$ be the convex hull of $(a_1 \vect{e}_1,\nv), \ldots, (a_{r-1}
  \vect{e}_{r-1},\nv)$, and $(a_r \vect{e}_r, a_r P)$, and let $Q_2$
  be given by $\{ (\vect{x},\vect{y}) \suchthat \vect{x} \in
  \sx(\vect{a}) \text{ and } x_r \eins + \vect{Ay} \ge 0 \}$.

  It is easy to see that $Q_1=Q_2$. The point $(\eins,\nv)$ is an
  interior point of $Q_2$, and all facets are at distance one. Since
  $Q_1$ has integral vertices, $Q_1=Q_2$ is reflexive.
\end{proof}
 From this the final result follows.
\begin{Corollary}
  Given a polytope $P$ the reflexive dimension
  of $kP$ is bounded from above by 
  $\refldim(kP) \le \refldim(P) + M \sqrt{\log k\,}.$
\end{Corollary}
\begin{proof}
From the proof of Proposition \ref{prop:upper}, we know that
the simplex $\sx(a_1,\ldots,a_{r-1}, k)$ is a reflexive
polytope of dimension $r$ bounded by $r \le M\sqrt{\log k\,}$, where
$M$ is a universal constant.
Let $Q$ be a reflexive polytope of dimension $\refldim(P)$ that
contains $P$ as a face.  Given $\sx$ and $Q$, it follows from
Proposition \ref{prop:polymult} that $kQ$ is a face of 
a reflexive polytope of dimension 
$\refldim(P) +r  \le \refldim(P) + M \sqrt{\log k\,}$.
\end{proof}
Let us conclude with a question.
Is 
the reflexive dimension of the Minkowski sum $P+P'$ 
bounded 
by $\refldim(P) + \refldim(P') + c$ ?
\subsection*{Acknowledgments}
This work was inspired by a question of Bernd Sturmfels who suggested
Definition~\ref{def:refldim}. We thank the OEIS~\cite{OEIS} for
pointing us to Egyptian fractions, \url{google} for pointing us to
David Eppstein, and David Eppstein for pointing us to~\cite{Vose}.

Finally, we thank Matthew Prior.


\bibliographystyle{plain}
\bibliography{alles,haasi}

\end{document}